\documentclass{article}
\usepackage{arXiv}

\hypersetup{
  pdftitle={Multivariate Lebesgue bounds for R-Leja lower grids},
  pdfauthor={Phil-Alexander Hofmann and Michael Hecht}
}

\usepackage{todonotes}

\usepackage{fancyhdr}
\title{\textbf{
    Multivariate Lebesgue bounds\\ for $\mathfrak R$-Leja lower grids
}}

\newcommand{\Nzero}{\mathbb N_0}
\newcommand{\oneA}{\mathbf 1_A}
\newcommand{\betaL}{\beta_d}

\author[1]{Phil-Alexander Hofmann\thanks{\href{mailto:p.hofmann@hzdr.de}{p.hofmann@hzdr.de}}}
\author[2]{Michael Hecht}

\affil[1]{Center for Advanced Systems Understanding,
Helmholtz-Zentrum Dresden-Rossendorf,
D-02826 G\"orlitz, Germany}
\affil[2]{Mathematical Institute, University of Wroc{\l}aw,
50-384 Wroc{\l}aw, Poland}

\begin{document}

\newpage

\maketitle
\thispagestyle{empty}
\begin{abstract}
We establish an algebraic upper bound for the Lebesgue constant of multivariate polynomial interpolation on $\mathfrak R$-Leja lower grids.
For each fixed dimension $d\geq2$ and every nonempty finite lower set $A\subset\mathbb N_0^d$, let $P_A$ denote the corresponding grid with $N=|A|$ nodes.
We prove $\Lambda(P_A)\leq C_dN^{\beta_d}$, where $q_d=d/(d-1)$ and $\beta_d=1+q_d^{-1}\log_2(1+2^{q_d})$.
The exponent satisfies $\beta_2=2.160964\ldots$ and increases to $\log_2 6=2.584962\ldots$ as $d\to\infty$, improving both the cubic estimate and the exponent $3-1/d$ obtained from quadratic prefix bounds.
Our proof exploits the binary structure of known univariate estimates to control sums of powers of Lebesgue constants over interpolation prefixes.
A weighted H\"older inequality and a disjoint-box packing argument then transfer these prefix-moment bounds to arbitrary lower sets.
The resulting estimate is uniform over all lower-set shapes and choices of disk Leja maximizers, with constants depending only on $d$.
\end{abstract}


\section{Introduction}
\label{sec:introduction}

We study the stability of multivariate polynomial interpolation on lower grids constructed via tensor products of $\mathfrak R$-Leja sequences indexed by finite lower sets.
A set $A\subset\mathbb{N}_0^d$ is a lower set if, for every $\bm{\nu}=(\nu_i)_{i=1,\ldots,d}\in A$, it also contains all $\bm{\mu}=(\mu_i)_{i=1,\ldots,d}$ with $0\leq\mu_i\leq\nu_i$ for $i=1,\ldots,d$.
The $\mathfrak R$-Leja sequences used here are obtained by taking the real parts of Leja sequences on the closed complex unit disk starting at $1$ and retaining distinct values in their order of first occurrence.
These sequences were introduced and studied by Calvi and Phung~\cite{CalviPhung2012}, building on their earlier work on Leja sequences for the unit disk and their applications to multivariate interpolation~\cite{CalviPhung2011}.

Our main result, Theorem~\ref{thm:main}, establishes a uniform algebraic upper bound for the Lebesgue constant of the resulting lower grids in terms of the total number $N=|A|$ of interpolation nodes, independently of the shape of the lower set.
For each fixed dimension $d\geq2$, the bound admits the exponent
\begin{equation}
    \label{eq:beta}
    \betaL\coloneqq 1+\frac1{q_d}\log_2(1+2^{q_d}),
    \qquad
    q_d\coloneqq\frac{d}{d-1}.
\end{equation}

\begin{theorem}[$\mathfrak R$-Leja grids on arbitrary lower sets]
    \label{thm:main}
    Let the dimension $d\geq2$ be fixed.
    Choose an $\mathfrak R$-Leja sequence $(p_{i,k})_{k\geq0}$ in each coordinate $i$, obtained from a Leja sequence on the complex unit disk starting at $1$.
    Further, for every nonempty finite lower set $A\subset\Nzero^d$, let
    \begin{equation}
        \label{eq:grid-space}
        P_A=\{(p_{1,\nu_1},\ldots,p_{d,\nu_d})\mid\bm{\nu}\in A\},
        \qquad
        \Pi_A=\operatorname{span}\{\bm{x}^{\bm{\nu}}\mid\bm{\nu}\in A\}.
    \end{equation}
    Let $\mathcal I_A$ be the interpolation operator from $C([-1,1]^d)$ into the downward closed polynomial space $\Pi_A$ at $P_A$.
    Writing $N\coloneqq|A|=|P_A|$ for the number of interpolation nodes, there exists a constant $C_d>0$ such that
    \begin{equation}
        \label{eq:main}
        \Lambda(P_A)=\|\mathcal I_A\|\leq C_dN^{\betaL},
    \end{equation}
    where $\|\mathcal I_A\|$ denotes the operator norm induced by the supremum norm on the cube $[-1,1]^d$.
    The constant depends only on $d$, and is independent of $A$ and the choices of disk Leja maximizers.
\end{theorem}

The Lebesgue constant quantifies the stability of polynomial interpolation and controls the interpolation error relative to the best approximation error.
For every $f\in C([-1,1]^d)$, the identity $\mathcal I_Ap=p$ for $p\in\Pi_A$ yields the inequality
\begin{equation}
    \label{eq:lebesgue-approximation}
    \begin{aligned}
        \|f-\mathcal I_Af\|_\infty
        &\leq(1+\Lambda(P_A))\inf_{p\in\Pi_A}\|f-p\|_\infty\\
        &\leq(1+C_dN^{\betaL})\inf_{p\in\Pi_A}\|f-p\|_\infty.
    \end{aligned}
\end{equation}
Thus, a smaller exponent reduces the polynomial factor relating the interpolation error to the best approximation error in $\Pi_A$.

\paragraph{Contribution and relation to previous bounds.}
Estimating the hierarchical increments separately gives $O(N^3)$ from quadratic coordinate-prefix bounds \cite[Lemma~3.1]{ChkifaCohenSchwab2014}.
Cohen and Chkifa give this cubic estimate for $\mathfrak R$-Leja grids uniformly over dimensions and lower-set shapes \cite[Corollary~5.5]{CohenChkifa2015}.
This can be refined to $O_d(N^{3-1/d})$ with a constant depending on the dimension $d$, as demonstrated in Section~\ref{sec:quadratic-comparison}, via collecting tensor prefixes and using the lower-set geometry.
For $d=2$, this specializes to $O(N^{5/2})$, which also follows from the standard staircase formula and quadratic prefix bounds \cite{CarnicerGodes2015,Chkifa2015Bounds}.
Theorem~\ref{thm:main} reduces the exponent further to $\beta_d$ via the prefix-moments of Lemma~\ref{lem:moment}, as summarized in Table~\ref{tab:exponents}.
In particular, for every fixed $d$, the exponent in \eqref{eq:main} satisfies $\betaL<\log_2 6<3$, while its constants may depend on $d$.
This bound is uniform, however, for particular hyperbolic shapes smaller bounds can be derived, as demonstrated in Section~\ref{sec:shape-bounds}.

\medskip
\begin{table}[!htbp]
    \centering
    \renewcommand{\arraystretch}{1.12}
    \begin{tabular}{@{}cccc@{}}
        \toprule
        Dimension & \shortstack{Hierarchical sum} & \shortstack{Geometric refinement} & \shortstack{Prefix moment refinement}\\
        $d$ & $3$ & $3-1/d$ & $\betaL$\\
        \midrule
        $2$ & $3$ & $2.50$ & $\mathbf{2.160964\ldots}$\\
        $3$ & $3$ & $2.\overline{6}$ & $\mathbf{2.291167\ldots}$\\
        \midrule
        $d\to\infty$ & $3$ & $3$ & $\log_2 6=2.584962\ldots$\\
        \bottomrule
    \end{tabular}
    \caption{
        Exponents $p$ in the bounds $\Lambda(P_A)=O_d(N^p)$, rounded to six decimal places.
        Here, $N$ is the total number of nodes.
        The last row gives the exponent limits as $d\to\infty$.
        The last two columns allow constants depending on $d$.
        All three comparisons are uniform over finite lower-set shapes.
    }
    \label{tab:exponents}
\end{table}

Chkifa's earlier analysis \cite[Theorem~3.1, equations~(3.11)--(3.12)]{Chkifa2013Lebesgue} establishes an $O(k^2\log k)$ bound for $\mathfrak R$-Leja prefixes and recalls the general hierarchical estimate.
The subsequent estimate \eqref{eq:source} \cite{Chkifa2015Bounds} captures the dependence on the binary expansion of the prefix length.
We exploit this structure to derive univariate prefix-moment bounds, which we then transfer to arbitrary finite lower sets using the tensor-prefix representation \cite{DynFloater2014}, a weighted product bound, and disjoint-box packing.

\paragraph{The proof in one line.}
The argument goes as follows, where $L_i(a)$ is the Lebesgue constant for interpolation at the first $a$ nodes in coordinate $i$, and $n$ is an arbitrary one-dimensional prefix-length cutoff,
\begin{equation}
    \label{eq:proof-roadmap}
    \text{Chkifa's binary prefix estimate}
    \ \Longrightarrow\
    \sum_{a=1}^nL_i(a)^{q_d}\lesssim n^{q_d\betaL}
    \ \Longrightarrow\
    \Lambda(P_A)\lesssim_d N^{\betaL}.
\end{equation}
The second implication specifically uses the lower-set geometry.

\paragraph{Organization of the paper.}
The remainder of this paper is organized as follows.
Section~\ref{sec:staircase} derives a tensor-prefix representation of the multivariate interpolation operator and discusses preliminary bounds of the Lebesgue constants.
Section~\ref{sec:prefix-moments} establishes the univariate prefix-moment bounds and Section~\ref{sec:moment-transfer} transfers them to arbitrary multivariate lower grids.
Finally, Section~\ref{sec:numerical-validation} compares the theoretical bounds with numerical estimates of the Lebesgue constants.

\section{Tensor-prefix representation and preliminary bounds}
\label{sec:staircase}

For each coordinate $i$ and integer $a\geq 1$, let $J_i(a)$ denote the interpolation operator at the first $a$ nodes $p_{i,0},\ldots,p_{i,a-1}$ onto polynomials of degree at most $a-1$.
Its Lebesgue constant is $L_i(a)=\|J_i(a)\|$, where the operator norm is induced by the supremum norm on $[-1,1]$.

Define the maximal prefix lengths by
\begin{equation}
    \label{eq:coordinate-prefix-lengths}
    n_i=1+\max_{\bm{\nu}\in A}\nu_i,\qquad 1\leq i\leq d,
    \qquad N\leq\prod_{i=1}^d n_i.
\end{equation}
The full tensor-product grid contains $\prod_{i=1}^d n_i$ nodes, whereas the lower grid associated with a nonempty finite lower set $A \subset \Nzero^d$ contains exactly $N=|A|$ nodes, while equality holds precisely when $A=\times_{i=1}^d\{0,\ldots,n_i-1\}$.
Distinct coordinate nodes and downward closure ensure that the lower grid is unisolvent, as follows from the triangular structure of its Newton basis~\cite{Hecht2026}.

\subsection{The tensor-prefix representation of the interpolation operator}
\label{sec:staircase-representation}

We express $\mathcal I_A$ as a linear combination of tensor-product interpolation operators on boxes in $A$.

Set $J_i(0)=0$ and define
\begin{equation}
    \label{eq:increments}
    \Delta_i(k)=J_i(k)-J_i(k-1),
    \qquad
    J_i(a)=\sum_{k=1}^{a}\Delta_i(k).
\end{equation}
The increment $\Delta_i(k)$ vanishes at the first $k-1$ nodes, since both interpolants coincide there.
The standard hierarchical representation is
\cite[equation~(2.16) and Theorem~2.1]{ChkifaCohenSchwab2014}
\begin{equation}
    \label{eq:hierarchical}
    \mathcal I_A
    =\sum_{\bm{\mu}\in A}\bigotimes_{i=1}^d\Delta_i(\mu_i+1).
\end{equation}

Expanding the differences in \eqref{eq:hierarchical} gives
\begin{equation}
    \label{eq:inclusion-exclusion}
    \mathcal I_A
    =
    \sum_{\bm{\varepsilon}\in\{0,1\}^d}
    (-1)^{|\bm{\varepsilon}|_1}
    \sum_{\bm{\mu}\in A}
    \bigotimes_{i=1}^d J_i(\mu_i+1-\varepsilon_i).
\end{equation}
Here $\varepsilon_i=0$ selects $J_i(\mu_i+1)$, while
$\varepsilon_i=1$ selects $-J_i(\mu_i)$.
For each fixed $\bm{\varepsilon}$, substitute
$\bm{\nu}=\bm{\mu}-\bm{\varepsilon}$ in the inner sum.
The tensor factors then become $J_i(\nu_i+1)$, and the restriction
$\bm{\mu}\in A$ becomes $\bm{\nu}+\bm{\varepsilon}\in A$.
Terms with some $\nu_i=-1$ vanish because $J_i(0)=0$.
Writing $\oneA$ for the indicator function of $A$, we therefore obtain
\begin{equation}
    \label{eq:inclusion-exclusion-shifted}
    \mathcal I_A
    =
    \sum_{\bm{\varepsilon}\in\{0,1\}^d}
    (-1)^{|\bm{\varepsilon}|_1}
    \sum_{\bm{\nu}\in\Nzero^d}
    \oneA(\bm{\nu}+\bm{\varepsilon})
    \bigotimes_{i=1}^d J_i(\nu_i+1).
\end{equation}
The tensor product no longer depends on $\bm{\varepsilon}$.
Interchanging the sums, which have only finitely many nonzero terms,
and collecting its coefficient yields
\begin{equation}
    \label{eq:combination}
    \begin{aligned}
        c_A(\bm{\nu})
        &=
        \sum_{\bm{\varepsilon}\in\{0,1\}^d}
        (-1)^{|\bm{\varepsilon}|_1}\oneA(\bm{\nu}+\bm{\varepsilon}),
        \\
        \mathcal I_A
        &=
        \sum_{\bm{\nu}\in\Nzero^d}
        c_A(\bm{\nu})\bigotimes_{i=1}^d J_i(\nu_i+1).
    \end{aligned}
\end{equation}
This is the tensor-prefix formula of Dyn and Floater \cite[Theorem~4, equations~(9)--(10)]{DynFloater2014}.

For a vector of positive integers $\bm{a}\in\{1,2,\ldots\}^d$, let the index box
\begin{equation}
    \label{eq:index-box}
    B(\bm{a})=\prod_{i=1}^d\{0,\ldots,a_i-1\}.
\end{equation}
If $\bm{\nu}\notin A$, downward closure excludes every $\bm{\nu}+\bm{\varepsilon}$, so $c_A(\bm{\nu})=0$.
Thus the tensor prefixes in \eqref{eq:combination} use only boxes contained in $A$:
\begin{equation}
    \label{eq:rectangle-below}
    \bm{\nu}\in A
    \quad\Longrightarrow\quad
    B(\bm{\nu}+\bm{1})\subset A,
    \qquad
    \prod_{i=1}^d(\nu_i+1)\leq N,
\end{equation}
where $\bm{1}=(1,\ldots,1)\in\Nzero^d$.
The coefficients $c_A(\bm{\nu})$ are integers, but their absolute values need not be equal to one.

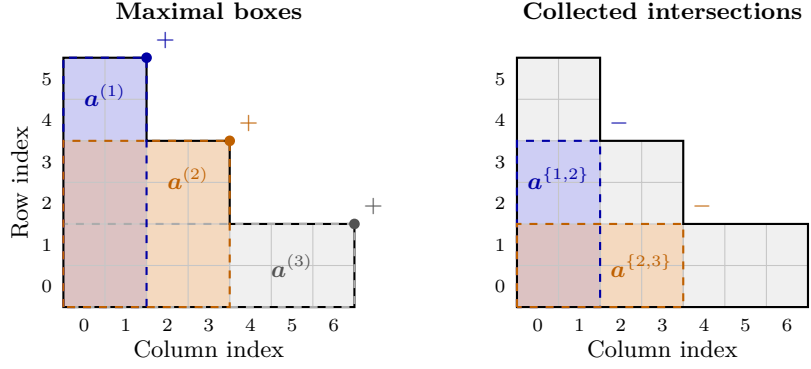
\begin{figure}[tbp]
    \centering
    \begin{tikzpicture}[x=0.55cm,y=0.55cm,font=\small]
        \begin{scope}
            \node[anchor=south] at (3.5,6.7){\textbf{Maximal boxes}};
            \fill[gray!12] (0,0)--(7,0)--(7,2)--(4,2)--(4,4)--(2,4)--(2,6)--(0,6)--cycle;
            \fill[blue!50,opacity=0.3] (0,0) rectangle (2,6);
            \fill[orange!70,opacity=0.3] (0,0) rectangle (4,4);
            \draw[step=1,gray!45,thin] (0,0) grid (2,6);
            \draw[step=1,gray!45,thin] (2,0) grid (4,4);
            \draw[step=1,gray!45,thin] (4,0) grid (7,2);
            \draw[thick] (0,0)--(7,0)--(7,2)--(4,2)--(4,4)--(2,4)--(2,6)--(0,6)--cycle;
            \draw[gray!65,thick,dashed] (0,0) rectangle (7,2);
            \draw[blue!65!black,thick,dashed] (0,0) rectangle (2,6);
            \draw[orange!75!black,thick,dashed] (0,0) rectangle (4,4);
            \node[blue!65!black] at (1,5.1){$\bm{a}^{(1)}$};
            \node[orange!75!black] at (3,3.1){$\bm{a}^{(2)}$};
            \node[gray!65!black] at (5.5,1){$\bm{a}^{(3)}$};
            \foreach \a/\b/\c in {2/6/blue!65!black,4/4/orange!75!black,7/2/gray!65!black}
            {\fill[\c] (\a,\b) circle (2pt);
            \node[above right,\c] at (\a,\b){$+$};}
            \foreach \k in {0,...,6}
                {\node[below,font=\scriptsize] at (\k+0.5,-0.05){$\k$};}
            \foreach \k in {0,...,5}
                {\node[left,font=\scriptsize] at (-0.05,\k+0.5){$\k$};}
            \node at (3.5,-1.0){Column index};
            \node[rotate=90] at (-1.05,3){Row index};
        \end{scope}
        \begin{scope}[xshift=6cm]
            \node[anchor=south] at (3.5,6.7){\textbf{Collected intersections}};
            \fill[gray!12] (0,0)--(7,0)--(7,2)--(4,2)--(4,4)--(2,4)--(2,6)--(0,6)--cycle;
            \fill[blue!50,opacity=0.3] (0,0) rectangle (2,4);
            \fill[orange!70,opacity=0.3] (0,0) rectangle (4,2);
            \draw[step=1,gray!45,thin] (0,0) grid (2,6);
            \draw[step=1,gray!45,thin] (2,0) grid (4,4);
            \draw[step=1,gray!45,thin] (4,0) grid (7,2);
            \draw[thick] (0,0)--(7,0)--(7,2)--(4,2)--(4,4)--(2,4)--(2,6)--(0,6)--cycle;
            \draw[blue!65!black,thick,dashed] (0,0) rectangle (2,4);
            \draw[orange!75!black,thick,dashed] (0,0) rectangle (4,2);
            \node[blue!65!black] at (1,3.1){$\bm{a}^{\{1,2\}}$};
            \node[orange!75!black] at (3,1){$\bm{a}^{\{2,3\}}$};
            \node[above right,blue!65!black] at (2,4){$-$};
            \node[above right,orange!75!black] at (4,2){$-$};
            \foreach \k in {0,...,6}
                {\node[below,font=\scriptsize] at (\k+0.5,-0.05){$\k$};}
            \foreach \k in {0,...,5}
                {\node[left,font=\scriptsize] at (-0.05,\k+0.5){$\k$};}
            \node at (3.5,-1.0){Column index};
        \end{scope}
    \end{tikzpicture}
    \caption{
        An index diagram illustrating \eqref{eq:staircase} for $N=26$.
        The maximal prefix vectors are $\bm{a}^{(1)}=(2,6)$, $\bm{a}^{(2)}=(4,4)$ and $\bm{a}^{(3)}=(7,2)$.
        After collecting intersections, the three boxes on the left have coefficient $+1$, and the boxes with prefix vectors $(2,4)$ and $(4,2)$ on the right have coefficient $-1$.
        Each unit square represents one index.
        Gray cells also belong to $A$.
    }
    \label{fig:staircase}
\end{figure}

For the geometric interpretation in Figure~\ref{fig:staircase}, let
$\bm{\mu}^{(1)},\ldots,\bm{\mu}^{(t)}\in\Nzero^d$ be the distinct maximal indices of $A$
with respect to the coordinatewise order, and set $\bm{a}^{(j)}=\bm{\mu}^{(j)}+\bm{1}$.
Every index in $A$ lies below a maximal index, hence
\begin{equation}
    \label{eq:maximal-boxes}
    A=\bigcup_{j=1}^t B(\bm{a}^{(j)}).
\end{equation}
For each nonempty $F\subseteq\{1,\ldots,t\}$, define $a_i^F=\min_{j\in F}a_i^{(j)}$.
Then
\begin{equation}
    \label{eq:intersection-boxes}
    B(\bm{a}^F)=\bigcap_{j\in F}B(\bm{a}^{(j)}).
\end{equation}
Applying inclusion--exclusion to the union of maximal boxes in \eqref{eq:maximal-boxes}, and telescoping the hierarchical sum over each intersection by \eqref{eq:increments}, gives the staircase representation
\begin{equation}
    \label{eq:staircase}
    \mathcal I_A
    =\sum_{\varnothing\ne F\subseteq\{1,\ldots,t\}}
        (-1)^{|F|+1}\mathcal I_{B(\bm{a}^F)},
    \qquad
    \mathcal I_{B(\bm{a})}=\bigotimes_{i=1}^d J_i(a_i).
\end{equation}
Here the singleton sets $F$ give the maximal boxes, and the intersection terms correct their overlaps.
Collecting equal intersections gives
\begin{equation}
    \label{eq:collected-box-coefficients}
    c_A(\bm{a}-\bm{1})
    =\sum_{\substack{\varnothing\ne F\subseteq\{1,\ldots,t\}\\\bm{a}^F=\bm{a}}}
        (-1)^{|F|+1}.
\end{equation}
This agrees with \eqref{eq:combination}. We write $\mathcal B_A$ for the distinct intersection boxes whose collected coefficient is nonzero.

\subsection{The cubic bound and its geometric refinement}
\label{sec:quadratic-comparison}

Chkifa's univariate estimate \cite{Chkifa2015Bounds} gives $L_i(a)\leq C_i a^2$ for interpolation at the first $a$ nodes of an $\mathfrak R$-Leja sequence, where $C_i>0$ is independent of $a$.
First, we bound the $N$ terms in \eqref{eq:hierarchical} separately.
With $L_i(0)=0$, the difference operators satisfy
\begin{equation}
    \label{eq:quadratic-increment}
    \|\Delta_i(k)\|
    \leq L_i(k)+L_i(k-1)
    \leq 2C_i k^2,\qquad k\geq1.
\end{equation}
Taking norms term by term in \eqref{eq:hierarchical} therefore gives
\begin{equation}
    \label{eq:cubic}
    \Lambda(P_A)
    \leq2^d\left(\prod_{i=1}^d C_i\right)\sum_{\bm{\nu}\in A} \left(\prod_{i=1}^d(\nu_i+1)\right)^2
    \leq2^d\left(\prod_{i=1}^d C_i\right)N^3.
\end{equation}
The box below each $\bm{\nu}\in A$ contains $\prod_{i=1}^d(\nu_i+1)\leq N$ indices.
Thus, each tensor increment has norm at most $2^d(\prod_{i=1}^d C_i)N^2$, and summing the $N$ terms gives the cubic bound.
This recovers the exponent $\theta+1$ of Chkifa, Cohen, and Schwab \cite[Lemma~3.1]{ChkifaCohenSchwab2014} for the quadratic case $\theta=2$.

Proposition~\ref{prop:transfer} improves the cubic bound by controlling the collected tensor prefixes through univariate moments.
Using only the quadratic prefix estimate and $q=q_d=d/(d-1)$ gives
\begin{equation}
    \label{eq:quadratic-moment}
    \left(\sum_{a=1}^nL_i(a)^q\right)^{1/q}
    \leq C_i\left(n\,n^{2q}\right)^{1/q}
    =C_i n^{3-1/d}.
\end{equation}
Applying the proposition with $K_i=C_i$ and $\gamma=3-1/d$ then yields
\begin{equation}
    \label{eq:quadratic-md}
    \Lambda(P_A)\leq C_d\left(\prod_{i=1}^dC_i\right)N^{3-1/d}.
\end{equation}
The geometric estimate therefore reduces the cubic exponent by $1/d$ in every fixed dimension, while the sharper prefix-moments developed in Section~\ref{sec:prefix-moments} further reduce the exponent to $\beta_d$.

\subsection{Bounds from the number and sizes of boxes}
\label{sec:shape-bounds}

The maximal-box representation \eqref{eq:staircase} gives a general bound after equal intersection boxes have been collected.
Let
\begin{equation}
    \label{eq:shape-box-count}
    M(A)=|\mathcal B_A|
\end{equation}
be the number of boxes in this collected representation.
This count includes maximal boxes and retained intersection boxes. For instance, in Figure~\ref{fig:staircase}, $M(A)=5$.
The coefficients are accounted for separately.
Formula~\eqref{eq:combination} has $2^{d-1}$ indicator terms of each sign, so $|c_A(\bm{\nu})|\leq2^{d-1}$.
For each box $B(\bm{a})$, the tensor-product Lebesgue function factorizes:
\begin{equation}
    \label{eq:tensor-term}
    \Lambda(P_{B(\bm{a})})
    =\|\mathcal I_{B(\bm{a})}\|
    =\prod_{i=1}^dL_i(a_i).
\end{equation}
Collecting equal intersections in \eqref{eq:staircase} and then taking norms gives the general bound
\begin{equation}
    \label{eq:shape-volume-sum}
    \Lambda(P_A)
    \leq\sum_{B(\bm{a})\in\mathcal B_A}|c_A(\bm{a}-\bm{1})|
        \Lambda(P_{B(\bm{a})}).
\end{equation}
This holds for every finite lower set and every choice of distinct nodes.
For the following two families of lower sets, we evaluate this sum using their box coefficients and the quadratic prefix estimates $L_i(a)\leq C_i a^2$.

\paragraph{Dyadic sparse sets.}
Consider
\begin{equation}
    \label{eq:shape-dyadic-set}
    A_n^{\mathrm{dyad}}
    =\bigcup_{\substack{\bm{\ell}\in\Nzero^d\\|\bm{\ell}|_1=n}}
        B(2^{\ell_1},\ldots,2^{\ell_d}),
    \qquad n\in\Nzero.
\end{equation}
Here $\ell_i$ specifies the prefix length $2^{\ell_i}$ in coordinate $i$, and $|\bm{\ell}|_1=\ell_1+\cdots+\ell_d$ is the sum of exponents.
The maximal boxes have $|\bm{\ell}|_1=n$ and each contain $2^n$ indices.
Grouping the hierarchical increments according to these powers of two and collecting tensor prefixes gives the Smolyak combination formula \cite[Remark~1(1)]{SickelUllrich2007},
\begin{equation}
    \label{eq:shape-dyadic-combination}
    \mathcal I_{A_n^{\mathrm{dyad}}}
    =\sum_{j=0}^{\min\{n,d-1\}}(-1)^j\binom{d-1}{j}
        \sum_{\substack{\bm{\ell}\in\Nzero^d\\|\bm{\ell}|_1=n-j}}
        \bigotimes_{i=1}^d J_i(2^{\ell_i}).
\end{equation}
Thus, the retained boxes have $|\bm{\ell}|_1=n-j$ for $0\leq j\leq\min\{n,d-1\}$.
For any integer $m\geq0$, the vectors $\bm{\ell}\in\Nzero^d$ with $|\bm{\ell}|_1=m$ correspond to strings of $m$ identical units split by $d-1$ separators into $d$ possibly empty groups.
Choosing the separator positions among the $m+d-1$ positions gives $\binom{m+d-1}{d-1}$ vectors.
For $d=2$, the vectors are $(k,m-k)$ with $k=0,\ldots,m$, giving $m+1$ choices.
Applying this count with $m=n-j$ gives
\begin{equation}
    \label{eq:shape-dyadic-box-count}
    M(A_n^{\mathrm{dyad}})
    =\sum_{j=0}^{\min\{n,d-1\}}\binom{n-j+d-1}{d-1}
    \asymp_d(n+1)^{d-1}.
\end{equation}
For a box with $|\bm{\ell}|_1=n-j$, the one-dimensional prefix estimates give
\begin{equation}
    \label{eq:shape-dyadic-box-norm}
    \prod_{i=1}^d L_i(2^{\ell_i})
    \leq\left(\prod_{i=1}^d C_i\right)4^{|\bm{\ell}|_1}
    =\left(\prod_{i=1}^d C_i\right)4^{n-j}.
\end{equation}
Inserting these box bounds and the coefficients in \eqref{eq:shape-dyadic-combination} into the general estimate \eqref{eq:shape-volume-sum} yields
\begin{equation}
    \label{eq:shape-dyadic-level-bound}
    \begin{aligned}
        \Lambda(P_{A_n^{\mathrm{dyad}}})
        &\leq\left(\prod_{i=1}^d C_i\right)
            \sum_{j=0}^{\min\{n,d-1\}}\binom{d-1}{j}
                \binom{n-j+d-1}{d-1}4^{n-j}\\
        &\lesssim_d4^n(n+1)^{d-1}.
    \end{aligned}
\end{equation}
The same rate follows from the Smolyak stability estimate \cite[Theorem~1(i), equation~(17)]{KangWilcox2017} using the quadratic univariate prefix bounds.
The distinct interpolation nodes are counted by \cite[Lemma~2(i)]{SickelUllrich2007}, and can also be found in \cite[Remark~2.4(1)]{MininiWasem2025}, namely
\begin{equation}
    \label{eq:shape-dyadic-cardinality}
    N=|A_n^{\mathrm{dyad}}|
    =\sum_{s=0}^{\min\{n,d-1\}}
        \binom{d-1}{s}\binom ns\,2^{n-s}
    \asymp_d2^n(n+1)^{d-1}.
\end{equation}
Together with \eqref{eq:shape-dyadic-box-count}, this gives
\begin{equation}
    \label{eq:shape-dyadic-box-count-nodes}
    M(A_n^{\mathrm{dyad}})\asymp_d(1+\log N)^{d-1}.
\end{equation}
Thus, for the dyadic sparse family, the box count grows as a fixed power of $1+\log N$.
Its Lebesgue bound becomes
\begin{equation}
    \label{eq:shape-dyadic-bound}
    \Lambda(P_{A_n^{\mathrm{dyad}}})
    \lesssim_d4^n(n+1)^{d-1}
    \asymp_d\frac{N^2}{(1+\log N)^{d-1}}.
\end{equation}
For every fixed $d$, this bound grows much slower than the uniform $N^{\beta_d}$ bound in Theorem~\ref{thm:main}.

Ultimately, this estimate can be even sharpened using the dyadic-prefix bound $L_i(2^\ell)\leq3\,2^\ell$ from \cite[equation~(4.10)]{Chkifa2015Bounds}.
Substituting this bound into \eqref{eq:shape-dyadic-combination} and using \eqref{eq:shape-dyadic-cardinality} gives
\begin{equation}
    \Lambda(P_{A_n^{\mathrm{dyad}}})
    \lesssim_d 2^n(n+1)^{d-1}
    \asymp_d N.
\end{equation}

\paragraph{Total-degree sets.}
For $A_n^{(1)}=\{\bm{\nu}\in\Nzero^d:|\bm{\nu}|_1\leq n\}$, adding the slack coordinate $\nu_{d+1}=n-|\bm{\nu}|_1$ gives a bijection with nonnegative integer vectors in $d+1$ coordinates whose entries sum to $n$.
The same separator count therefore gives $N=\binom{n+d}{d}\asymp_d(n+1)^d$.
The maximal boxes satisfy $|\bm{\nu}|_1=n$, and the boxes retained after collecting intersections satisfy $\max\{0,n-d+1\}\leq|\bm{\nu}|_1\leq n$.
Indeed, writing $j=n-|\bm{\nu}|_1$, formula~\eqref{eq:combination} gives
\begin{equation}
    \label{eq:shape-total-degree-coefficients}
    c_{A_n^{(1)}}(\bm{\nu})
    =\sum_{s=0}^j(-1)^s\binom ds
    =(-1)^j\binom{d-1}{j},
    \qquad 0\leq j\leq\min\{n,d-1\},
\end{equation}
and all other coefficients vanish.
For each $j$, there are $\binom{n-j+d-1}{d-1}$ vectors with $|\bm{\nu}|_1=n-j$, so
\begin{equation}
    \label{eq:shape-total-degree-box-count}
    M(A_n^{(1)})
    =\sum_{j=0}^{\min\{n,d-1\}}\binom{n-j+d-1}{d-1}
    \asymp_d(n+1)^{d-1}
    \asymp_d N^{1-1/d}.
\end{equation}
The total-degree family therefore has polynomially many retained boxes in terms of the node count.
For each retained box, \eqref{eq:tensor-term} and the quadratic prefix estimates give
$\Lambda(P_{B(\bm{\nu}+\bm{1})})\leq\left(\prod_{i=1}^d C_i\right)\prod_{i=1}^d(\nu_i+1)^2$.
Summing these box bounds in \eqref{eq:shape-volume-sum} yields
\begin{equation}
    \label{eq:shape-total-degree-bound}
    \begin{aligned}
        \Lambda(P_{A_n^{(1)}})
        &\leq\left(\prod_{i=1}^d C_i\right)
            \sum_{j=0}^{\min\{n,d-1\}}\binom{d-1}{j}
            \sum_{\substack{\bm{\nu}\in\Nzero^d\\|\bm{\nu}|_1=n-j}}
                \prod_{i=1}^d(\nu_i+1)^2\\
        &\lesssim_d(n+1)^{d-1}(n+1)^{2d}
            \asymp_d N^{3-1/d}.
    \end{aligned}
\end{equation}
Here the two powers of $n+1$ come from the box count and the bound $\nu_i+1\leq n+1$ for each coordinate, respectively.
For total-degree sets, the resulting exponent $3-1/d$ reproduces the general estimate \eqref{eq:quadratic-md} and remains larger than the exponent in Theorem~\ref{thm:main}.

These examples show how the box estimate benefits from sparse structure.
For denser families such as the total-degree sets, the number of retained boxes can make the sum of individual quadratic prefix bounds too coarse to improve the general estimate.
The usefulness of this approach therefore depends on both the number and the sizes of the boxes.
Next, to obtain a better bound for arbitrary lower sets, we estimate in Section~\ref{sec:prefix-moments}  the univariate prefix-moments, and combine these in Section~\ref{sec:moment-transfer} with the binary coefficient structure and geometric packing.

\section{Prefix-moments for $\mathfrak R$-Leja sequences}
\label{sec:prefix-moments}

To improve on the preceding bounds, we estimate sums of powers of the univariate Lebesgue constants over all prefixes up to a given length.
Lemma~\ref{lem:moment} establishes this prefix-moment bound for every exponent $q>0$. 
In Section~\ref{sec:moment-transfer}, we combine this bound with the geometry of lower sets to bound the multivariate Lebesgue constant in terms of the total node count $N$.
This transfer relies on the weighted product bound in Lemma~\ref{lem:holder}, choosing the specific exponent $q_d=d/(d-1)$.

Fix an $\mathfrak R$-Leja sequence $(p_k)_{k\geq0}$ in $[-1,1]$ and suppress the coordinate subscript $i$.
Its prefix of length $a$ consists of the first $a$ nodes $p_0,\ldots,p_{a-1}$.
Denote interpolation at these nodes by polynomials of degree at most $a-1$ by $J(a)$.
The Lebesgue constant is then given by
\begin{equation}
    \label{eq:univariate-lebesgue}
    L(a)
    =\|J(a)\|
    =\max_{x\in[-1,1]}\sum_{k=0}^{a-1}|\ell_{k,a}(x)|,
    \qquad
    a\geq1,
\end{equation}
where $\ell_{k,a}$ is the Lagrange polynomial which is equal to one at $p_k$ and zero at the other prefix nodes.
For $q>0$, we seek a bound for $\sum_{a=1}^nL(a)^q$, the prefix-moment of order $q$.

Chkifa organizes the prefix lengths into blocks between the lengths $2^r+1$ and $2^{r+1}+1$.
For an intermediate length, write
\begin{equation}
    \label{eq:prefix-offset}
    a=2^r+1+l,\qquad r\geq1,\qquad 1\leq l<2^r.
\end{equation}
Here $l$ counts the nodes added after the prefix of length $2^r+1$.
The estimate depends on two quantities of its binary expansion.
The digit sum $\sigma_1(l)$ is the number of digits equal to one, with $\sigma_1(0)=0$.
The integer $\sigma_0(l)$ is the number of zeros at the right end, equivalently the largest $\sigma$ such that $2^\sigma$ divides $l$.
For example,
\begin{equation}
    \label{eq:binary-example}
    6=(110)_2,\qquad \sigma_1(6)=2,\qquad \sigma_0(6)=1.
\end{equation}

The estimate uses an auxiliary Lebesgue constant $U_j$ for interpolation at the $j$ equally spaced roots of unity on the complex unit circle, with the maximum norm taken on the closed unit disk.
With this notation, Chkifa proved \cite[Theorem~3.2, equation~(3.12)]{Chkifa2015Bounds}
\begin{equation}
    \label{eq:source}
    L(2^r+1+l)\leq6\sqrt5\,
     2^{r+\sigma_1(l)-\sigma_0(l)}U_{2^{\sigma_0(l)}},\qquad 1\leq l<2^r.
\end{equation}
To bound $U$, let $\zeta_h=e^{2\pi\mathrm{i}h/j}$, for $0\leq h<j$.
For these roots, let $u_{h,j}$ denote the $h$-th Lagrange basis polynomial, i.e., the unique polynomial of degree at most $j-1$ satisfying $u_{h,j}(\zeta_k)=\delta_{hk}$.
For equally spaced roots of unity, these polynomials admit the explicit representation
\begin{equation}
    \label{eq:roots-cardinal}
    u_{h,j}(z) = \frac1j\sum_{t=0}^{j-1} \left(\frac{z}{\zeta_h}\right)^t,
    \qquad
    u_{h,j}(\zeta_k) = \delta_{hk}.
\end{equation}
Indeed, if $k=h$, all $j$ summands equal one, giving $u_{h,j}(\zeta_h)=j/j=1$.
If $k\ne h$, set $\rho=\zeta_k/\zeta_h$.
Then $\rho\ne1$, but $\rho^j=1$, so the geometric sum gives
\begin{equation}
    \label{eq:roots-cardinal-values}
    u_{h,j}(\zeta_k)
    = \frac1j\sum_{t=0}^{j-1}\rho^t
    = \frac1j\frac{1-\rho^j}{1-\rho}=0.
\end{equation}
Consequently, $\sum_{h=0}^{j-1}f(\zeta_h)u_{h,j}(z)$ interpolates the values of $f$ at $\zeta_h$.
Since $|\zeta_h|=1$, the triangle inequality further gives
\begin{equation}
    \label{eq:roots-cardinal-bound}
    |u_{h,j}(z)| \leq \frac1j\sum_{t=0}^{j-1}|z|^t\leq1,
    \qquad
    |z| \leq 1.
\end{equation}
Hence, the Lebesgue function sums the absolute values of all $j$ Lagrange polynomials, delivering
\begin{equation}
    \label{eq:roots-lebesgue-bound}
    U_j = \max_{|z|\leq1}\sum_{h=0}^{j-1}|u_{h,j}(z)|\leq j.
\end{equation}
In particular, $U_{2^{\sigma_0(l)}}\leq2^{\sigma_0(l)}$ cancels the factor $2^{-\sigma_0(l)}$ in \eqref{eq:source}, yielding
\begin{equation}
    \label{eq:envelope}
    L(2^r+1+l) \leq 6\sqrt5\,2^r2^{\sigma_1(l)}.
\end{equation}
The factor $2^{\sigma_1(l)}$ attains its maximum $2^r$ only when $l=2^r-1$.


\begin{lemma}[Univariate prefix-moment bound]
    \label{lem:moment}
    For every $q>0$, there is a constant $K_q>0$, uniform over all $\mathfrak R$-Leja sequences, such that for all integers $n\geq1$,
    \begin{equation}
        \label{eq:moment}
        \sum_{a=1}^n L(a)^q\leq K_q n^{\alpha_q},
        \qquad \alpha_q\coloneqq q+\log_2(1+2^q).
    \end{equation}
\end{lemma}

\begin{proof}
    Raising \eqref{eq:envelope} to the power $q$ gives
    \begin{equation}
        \label{eq:squared-envelope}
        L(2^r+1+l)^q
        \leq (6\sqrt5)^q(2^r)^q\bigl(2^{\sigma_1(l)}\bigr)^q
        = (6\sqrt5)^q2^{qr}2^{q\sigma_1(l)},
        \qquad
        1\leq l<2^r.
    \end{equation}
    To sum these powers of prefix bounds, write
    \begin{equation}
        \label{eq:binary-digits}
        l = \sum_{j=0}^{r-1}\varepsilon_j2^j,
        \qquad \varepsilon_j\in\{0,1\},
        \qquad \sigma_1(l)
        = \sum_{j=0}^{r-1}\varepsilon_j.
    \end{equation}
    As $l$ runs from $0$ to $2^r-1$, its $r$ binary digits run through every possible choice, including leading zeros.
    Since $2^{q\sigma_1(l)} = \prod_{j=0}^{r-1}2^{q\varepsilon_j}$, we may rewrite the sum as
    \begin{equation}
        \label{eq:binary-reindex}
        \sum_{l=0}^{2^r-1}2^{q\sigma_1(l)}
        = \sum_{(\varepsilon_0,\ldots,\varepsilon_{r-1})\in\{0,1\}^r} \prod_{j=0}^{r-1}2^{q\varepsilon_j}.
    \end{equation}
    To evaluate this sum, first keep $\varepsilon_1,\ldots, \varepsilon_{r-1}$ fixed and sum only over $\varepsilon_0=0,1$.
    The two terms have the same remaining factor, so
    \begin{equation}
        \label{eq:binary-one-digit}
        \sum_{\varepsilon_0=0}^1\prod_{j=0}^{r-1}2^{q\varepsilon_j}
        = 1\cdot\prod_{j=1}^{r-1}2^{q\varepsilon_j} +2^q\cdot\prod_{j=1}^{r-1}2^{q\varepsilon_j}
        = (1+2^q)\prod_{j=1}^{r-1}2^{q\varepsilon_j}.
    \end{equation}
    Summing over each remaining digit in the same way introduces one factor $1+2^q$ per digit.
    There are $r$ digits, hence
    \begin{equation}
        \label{eq:binary}
        \sum_{l=0}^{2^r-1}2^{q\sigma_1(l)}
        = \prod_{j=0}^{r-1}((2^q)^0+(2^q)^1)
        = (1+2^q)^r.
    \end{equation}
    These weighted digit sums also appear in \cite[Appendix~A, equation~(A.3)]{Chkifa2022Newton}, with $b=2^q$.
    Summing \eqref{eq:squared-envelope} over the block and adding the nonnegative $l=0$ term on the right, \eqref{eq:binary} gives
    \begin{equation}
        \label{eq:block}
        \sum_{l=1}^{2^r-1}L(2^r+1+l)^q
        \leq (6\sqrt5)^q2^{qr}\sum_{l=0}^{2^r-1}2^{q\sigma_1(l)}
        = (6\sqrt5)^q2^{\alpha_qr}.
    \end{equation}
    The sum in \eqref{eq:block} starts at $l=1$, so it omits the complete prefix of length $2^r+1$.
    Include that missing term by grouping the prefix lengths into the full blocks
    \begin{equation}
        \label{eq:full-prefix-blocks}
        \mathcal B_r
        = \{2^r+1,\ldots,2^{r+1}\}
        = \{2^r+1+l\mid0\leq l<2^r\},\qquad r\geq0.
    \end{equation}
    Thus $\mathcal B_0 = \{2\}$, $\mathcal B_1 = \{3,4\}$, $\mathcal B_2 = \{5,6,7,8\}$, and so on.
    These consecutive, disjoint blocks contain every prefix length $a\geq2$.

    For the first term of each block, \cite[equation~(3.8)]{Chkifa2015Bounds} gives $L(2^r+1)\leq C_0(r+1)$ with an absolute constant $C_0>0$, enlarged if necessary to include $r=0$.
    Since $r+1\leq2^r$ for all integers $r\geq0$,
    \begin{equation}
        \label{eq:complete-prefix-square}
        L(2^r+1)^q
        \leq C_0^q(r+1)^q
        \leq C_0^q2^{qr}
        \leq C_0^q2^{\alpha_qr}.
    \end{equation}
    Adding this term to \eqref{eq:block} gives the full-block bound
    \begin{equation}
        \label{eq:full-prefix-block-bound}
        \sum_{a\in\mathcal B_r}L(a)^q
        = L(2^r+1)^q+\sum_{l=1}^{2^r-1}L(2^r+1+l)^q
        \leq C_q2^{\alpha_qr},
        \qquad
        C_q = C_0^q+(6\sqrt5)^q.
    \end{equation}
    For $r=0$, the second sum is empty, so this bound also applies to $\mathcal B_0$.

    It remains to sum these block bounds up to an arbitrary length $n$.
    For $n\geq2$, let $R\geq0$ be the unique integer with
    \begin{equation}
        \label{eq:final-prefix-block}
        2^R+1 \leq n \leq 2^{R+1}.
    \end{equation}
    Then $n$ lies in $\mathcal B_R$.
    All earlier blocks are included in full, while only the lengths up to $n$ in the last block are needed.
    Since the summands are nonnegative, enlarging that last sum to the whole block only increases it. 
    Using $L(1)=1$ and \eqref{eq:full-prefix-block-bound},
    \begin{align}
        \label{eq:prefix-block-sum}
        \sum_{a=1}^n L(a)^q
        & = 1+\sum_{r=0}^{R-1}\sum_{a\in\mathcal B_r}L(a)^q +\sum_{\substack{a\in\mathcal B_R\\a\leq n}}L(a)^q\\
        \label{eq:prefix-block-majorant}
        &\leq 1+C_q\sum_{r=0}^R2^{\alpha_qr}.
    \end{align}
    The sum over earlier blocks is empty when $R=0$.
    Evaluating the geometric sum gives
    \begin{equation}
        \label{eq:prefix-geometric-sum}
        1+C_q\sum_{r=0}^R2^{\alpha_qr}
        = 1+C_q\frac{2^{\alpha_q(R+1)}-1}{2^{\alpha_q}-1}
        \leq K_q2^{\alpha_qR},
        \qquad
        K_q = 1+\frac{2^{\alpha_q}C_q}{2^{\alpha_q}-1}.
    \end{equation}
    Finally, \eqref{eq:final-prefix-block} implies $2^R<n$, hence
    \begin{equation}
        \label{eq:prefix-scale-to-length}
        2^{\alpha_qR} = (2^R)^{\alpha_q}\leq n^{\alpha_q}.
    \end{equation}
    Combining these inequalities proves \eqref{eq:moment} for $n\geq2$.
    The case $n=1$ follows from $L(1)=1$ and $K_q\geq1$.
\end{proof}

For comparison, replacing every prefix by the coarse bound $L(a)\lesssim a^2$ gives only $\sum_{a=1}^nL(a)^q\lesssim n^{2q+1}$.
The detailed binary estimate reduces this moment exponent to $\alpha_q$.
Next, Section~\ref{sec:moment-transfer} transfers $\alpha_{q_d}/q_d=\betaL$ to the multivariate node count.

\section{Transferring the prefix-moments to arbitrary lower grids}
\label{sec:moment-transfer}

Starting from \eqref{eq:combination}, Lemma~\ref{lem:line} bounds the total absolute coefficient weight along each coordinate-parallel line. 
Lemma~\ref{lem:holder} establishes a general weighted product bound under such line-weight bounds. 
Proposition~\ref{prop:transfer} then combines this estimate with the univariate prefix-moment bounds and a disjoint-box packing argument to obtain a bound for $\Lambda(P_A)$ in terms of the total node count $N=|A|$, preserving the exponent of the $q$th root of the univariate prefix-moment.

The coefficient estimates and the packing argument depend only on the lower set. 
Consequently, Proposition~\ref{prop:transfer} applies to arbitrary coordinate sequences of pairwise distinct nodes in $[-1,1]$, provided their Lebesgue constants satisfy the required univariate prefix-moment bound. 
For $\mathfrak R$-Leja sequences, this condition follows from Lemma~\ref{lem:moment}.

\subsection{The total coefficient weight on a line}
\label{sec:line-weight}

In \eqref{eq:combination}, all contributions with the same tensor prefix have already been combined into the coefficient $c_A(\bm{\nu})$.
Set $S_A=\{\bm{\nu}\in\Nzero^d\mid c_A(\bm{\nu})\ne0\}$ and $w_A(\bm{\nu})=|c_A(\bm{\nu})|$.
Applying the triangle inequality to that sum and using the tensor-product factorization in \eqref{eq:tensor-term} gives
\begin{equation}
    \label{eq:absolute-cost}
    \begin{aligned}
        \Lambda(P_A)=\|\mathcal I_A\|
        &\leq\sum_{\bm{\nu}\in S_A}|c_A(\bm{\nu})|
            \left\|\bigotimes_{i=1}^dJ_i(\nu_i+1)\right\|\\
        &=\sum_{\bm{\nu}\in S_A}w_A(\bm{\nu})\prod_{i=1}^dL_i(\nu_i+1).
    \end{aligned}
\end{equation}

Now, we bound the total coefficient weight $w_A$ along each coordinate-parallel line.

\begin{lemma}[Total coefficient weight on a line]
    \label{lem:line}
    Let $A\subset\Nzero^d$ be a nonempty finite lower set.
    Then, the support $S_A$ is contained in $A$.
    Further, fix $i\in\{1,\ldots,d\}$.
    For $\bm{\nu}\in\Nzero^d$, define the vector obtained by deleting coordinate $i$ by
    \begin{equation}
        \label{eq:coordinate-deletion}
        \bm{\nu}_{-i}\coloneqq
        (\nu_1,\ldots,\nu_{i-1},\nu_{i+1},\ldots,\nu_d)
        \in\Nzero^{d-1}.
    \end{equation}
    For $\bm{\eta}\in\Nzero^{d-1}$ and $t\in\Nzero$, define the insertion map
    \begin{equation}
        \label{eq:coordinate-insertion}
        \Phi_i(\bm{\eta},t)\coloneqq
        (\eta_1,\ldots,\eta_{i-1},t,\eta_i,\ldots,\eta_{d-1})
        \in\Nzero^d.
    \end{equation}
    Thus, $\Phi_i(\bm{\nu}_{-i},\nu_i)=\bm{\nu}$.
    Then, for every fixed $\bm{\eta}\in\Nzero^{d-1}$, it holds
    \begin{equation}
        \label{eq:line-weight}
        \sum_{t\in\Nzero}
         w_A\bigl(\Phi_i(\bm{\eta},t)\bigr)
        \leq 2^{d-1}.
    \end{equation}
\end{lemma}

\begin{proof}
    If $\bm{\nu}\notin A$, downward closure excludes every $\bm{\nu}+\bm{\varepsilon}$, so $c_A(\bm{\nu})=0$.
    For the line bound, fix $\bm{\eta}=\bm{\nu}_{-i}$ and split the sum in \eqref{eq:combination} according to $b=\varepsilon_i\in\{0,1\}$ and $\bm{\varepsilon}'=\bm{\varepsilon}_{-i}\in\{0,1\}^{d-1}$.
    Then $\bm{\varepsilon}=\Phi_i(\bm{\varepsilon}',b)$, so $|\bm{\varepsilon}|_1=|\bm{\varepsilon}'|_1+b$ and
    $\Phi_i(\bm{\eta},t)+\bm{\varepsilon}=\Phi_i(\bm{\eta}+\bm{\varepsilon}',t+b)$.
    Substitution into \eqref{eq:combination} and evaluation of the inner sum over $b$ give
    \begin{equation}
        \label{eq:line-difference}
        \begin{aligned}
            c_A\bigl(\Phi_i(\bm{\eta},t)\bigr)
            &=\sum_{\bm{\varepsilon}'\in\{0,1\}^{d-1}}
                \sum_{b=0}^1(-1)^{|\bm{\varepsilon}'|_1+b}
                \oneA\bigl(\Phi_i(\bm{\eta}+\bm{\varepsilon}',t+b)\bigr)\\
            &=\sum_{\bm{\varepsilon}'\in\{0,1\}^{d-1}}(-1)^{|\bm{\varepsilon}'|_1}
                \bigl[\oneA\bigl(\Phi_i(\bm{\eta}+\bm{\varepsilon}',t)\bigr)
                -\oneA\bigl(\Phi_i(\bm{\eta}+\bm{\varepsilon}',t+1)\bigr)\bigr].
        \end{aligned}
    \end{equation}
    The two terms in brackets correspond to $b=0$ and $b=1$.
    For each fixed $\bm{\varepsilon}'$, downward closure makes the section
    $\{t\in\Nzero\mid\Phi_i(\bm{\eta}+\bm{\varepsilon}',t)\in A\}$ empty or an initial interval $\{0,\ldots,H-1\}$.
    The bracketed difference is therefore nonnegative and is one only at $t=H-1$ for a nonempty section. It is zero everywhere for an empty section.
    Its sum over $t\in\Nzero$ is at most one.
    Applying the triangle inequality to \eqref{eq:line-difference} and summing over the $2^{d-1}$ choices of $\bm{\varepsilon}'$ proves the assertion.
\end{proof}

As we will see next, this line weight bound provides the key ingredient for a weighted product bound, which allows us to control the multivariate sum in \eqref{eq:absolute-cost} through univariate prefix moments.

\subsection{A weighted product bound}
\label{sec:weighted-product}

We combine the line-weight bound with H\"older's inequality to control the weighted products in \eqref{eq:absolute-cost} by univariate moment sums.
This is the product estimate used later in the proof of Proposition~\ref{prop:transfer}.

\begin{lemma}[Weighted product bound]
    \label{lem:holder}
    Let $d\geq2$, let $w$ be a nonnegative function of finite support on $\Nzero^d$, and let $B\geq0$.
    Suppose that for every $i\in\{1,\ldots,d\}$ and every $\bm{\eta}\in\Nzero^{d-1}$, we have
    \begin{equation}
        \label{eq:holder-line-assumption}
        \sum_{t\in\Nzero}w\bigl(\Phi_i(\bm{\eta},t)\bigr)\leq B.
    \end{equation}
    Set $q=d/(d-1)$.
    Then, for nonnegative finitely supported sequences $f_i$ on $\Nzero$, $1\leq i\leq d$, it holds
    \begin{equation}
        \label{eq:weighted-holder}
        \sum_{\bm{\nu}\in\Nzero^d} w(\bm{\nu})\prod_{i=1}^df_i(\nu_i)
        \leq B\prod_{i=1}^d
                   \left(\sum_{k\in\Nzero}f_i(k)^q\right)^{1/q}.
    \end{equation}
\end{lemma}

\begin{proof}
    For $\bm{\nu}\in\Nzero^d$ and $1\leq i\leq d$, set
    \begin{equation}
        \label{eq:holder-factors}
        h_i(\bm{\nu})=w(\bm{\nu})^{1/d}
            \prod_{\substack{j=1\\j\ne i}}^d f_j(\nu_j)^{1/(d-1)}.
    \end{equation}
    Each $f_j(\nu_j)$ occurs in $d-1$ of these factors, so
    \begin{equation}
        \label{eq:holder-factor-identities}
        \prod_{i=1}^d h_i(\bm{\nu})
        =w(\bm{\nu})\prod_{i=1}^d f_i(\nu_i),
        \qquad
        h_i(\bm{\nu})^d
        =w(\bm{\nu})\prod_{\substack{j=1\\j\ne i}}^d f_j(\nu_j)^q.
    \end{equation}
    Applying H\"older's inequality to the $d$ factors $h_i$ gives
    \begin{equation}
        \label{eq:holder-step}
        \begin{aligned}
            \sum_{\bm{\nu}\in\Nzero^d} w(\bm{\nu})\prod_{i=1}^d f_i(\nu_i)
            &=\sum_{\bm{\nu}\in\Nzero^d}\prod_{i=1}^d h_i(\bm{\nu})\\
            &\leq\prod_{i=1}^d
                \left(\sum_{\bm{\nu}\in\Nzero^d}h_i(\bm{\nu})^d\right)^{1/d}.
        \end{aligned}
    \end{equation}
    For each $i$, the factors $f_j(\nu_j)^q$ with $j\ne i$ are independent of $\nu_i$.
    Summing first along coordinate $i$ therefore gives
    \begin{equation}
        \label{eq:holder-fibre}
        \begin{aligned}
            \sum_{\bm{\nu}\in\Nzero^d}h_i(\bm{\nu})^d
            &=\sum_{\substack{\nu_j\in\Nzero\\1\leq j\leq d,\ j\ne i}}
                \left(\sum_{t\in\Nzero}w\bigl(\Phi_i(\bm{\nu}_{-i},t)\bigr)\right)
                \prod_{\substack{j=1\\j\ne i}}^d f_j(\nu_j)^q\\
            &\leq B\prod_{\substack{j=1\\j\ne i}}^d
                \left(\sum_{k\in\Nzero}f_j(k)^q\right).
        \end{aligned}
    \end{equation}
    The outer sum runs independently over all coordinates except $i$, so after applying the line bound the remaining sum factorizes.
    Substituting \eqref{eq:holder-fibre} into \eqref{eq:holder-step} gives $B^{d/d}=B$, and each univariate sum appears $d-1$ times with exponent $1/d$, proving \eqref{eq:weighted-holder} since $(d-1)/d=1/q$.
\end{proof}

For $d=2$, this is a weighted Cauchy--Schwarz inequality.
Then $q=2$, so the right-hand side of \eqref{eq:weighted-holder} becomes $B\left(\sum_{k\in\Nzero}f_1(k)^2\right)^{1/2}\left(\sum_{k\in\Nzero}f_2(k)^2\right)^{1/2}$.

\subsection{Transfer by disjoint interior boxes}
\label{sec:disjoint-box-transfer}

Finally, we combine the weighted product bound with a dyadic decomposition of the prefix lengths and a packing bound inside the lower set.

\begin{proposition}[Transfer from univariate prefix-moments]
    \label{prop:transfer}
    Let $d \in \mathbb N$ with $d \geq 2$, set $q=d/(d-1)$, and suppose $\gamma\geq1$.
    For each $i\in\{1,\ldots,d\}$, let $(p_{i,k})_{k\in\Nzero}$ be a sequence of pairwise distinct nodes in $[-1,1]$, and let $L_i(a)$ be the Lebesgue constant for interpolation at its first $a$ nodes.
    Assume that there are constants $K_i>0$, $1\leq i\leq d$, such that for every integer $n\geq1$, we have
    \begin{equation}
        \label{eq:moment-assumption}
        \left(\sum_{a=1}^nL_i(a)^q\right)^{1/q}
        \leq K_i n^\gamma,\qquad n\geq1,\quad 1\leq i\leq d.
    \end{equation}
    For every nonempty finite lower set $A\subset\Nzero^d$, let $P_A$ be the grid in \eqref{eq:grid-space} and set $N=|A|=|P_A|$.
    Then
    \begin{equation}
        \label{eq:transfer}
        \Lambda(P_A)\leq
        2^{(d-1)(1+2\gamma)}
              \left(\prod_{i=1}^dK_i\right)N^\gamma.
    \end{equation}
\end{proposition}

\begin{proof}
    We group the first $d-1$ prefix lengths dyadically and retain the exact height in the last coordinate to avoid an additional dyadic cutoff and packing loss.
    Any coordinate could serve as this height direction.
    Set
    \begin{equation}
        \label{eq:dyadic-prefix-bins}
        E_r=\{2^r,\ldots,2^{r+1}-1\},\qquad r\geq0.
    \end{equation}
    For $\bm{r}=(r_1,\ldots,r_{d-1})\in\Nzero^{d-1}$, let
    \begin{equation}
        \label{eq:partial-bins}
        S_{\bm{r}}=\{\bm{\nu}\in S_A\mid
                      \nu_i+1\in E_{r_i}\ \text{for }i<d\},\qquad
        R=r_1+\cdots+r_{d-1}.
    \end{equation}
    The nonempty groups are indexed by
    \begin{equation}
        \label{eq:nonempty-dyadic-groups}
        \mathcal R_A=\{\bm{r}\in\Nzero^{d-1}\mid S_{\bm{r}}\ne\varnothing\}.
    \end{equation}
    For each $\bm{r}\in\mathcal R_A$, define
    \begin{equation}
        \label{eq:bin-height}
        h_{\bm{r}}=\max_{\bm{\nu}\in S_{\bm{r}}}(\nu_d+1),
        \qquad v_{\bm{r}}=2^R h_{\bm{r}}.
    \end{equation}
    In particular, $h_{\bm{r}}\geq1$.

    \paragraph{Bound the whole group.}
    Use $f_i(k)=L_i(k+1)$ on the prefix range and zero elsewhere.
    For $1\leq i\leq d-1$, nonnegativity allows us to enlarge the dyadic block to a prefix starting at $a=1$.
    Applying \eqref{eq:moment-assumption} with $n=2^{r_i+1}$ gives
    \begin{equation}
        \label{eq:dyadic-block-moment}
        \left(\sum_{a\in E_{r_i}}L_i(a)^q\right)^{1/q}
        \leq\left(\sum_{a=1}^{2^{r_i+1}}L_i(a)^q\right)^{1/q}
        \leq K_i2^{\gamma(r_i+1)}.
    \end{equation}
    The last coordinate already uses a prefix, so \eqref{eq:moment-assumption} gives
    $\left(\sum_{a=1}^{h_{\bm{r}}}L_d(a)^q\right)^{1/q}\leq K_dh_{\bm{r}}^\gamma$.
    Applying Lemma~\ref{lem:holder} to $w_A$ restricted to $S_{\bm{r}}$, with the line-weight bound from Lemma~\ref{lem:line}, yields
    \begin{align}
        \label{eq:group-holder}
        \sum_{\bm{\nu}\in S_{\bm{r}}}
            w_A(\bm{\nu})\prod_{i=1}^dL_i(\nu_i+1)
        &\leq 2^{d-1}
              \prod_{i=1}^{d-1}\left(\sum_{a\in E_{r_i}}L_i(a)^q\right)^{1/q}
              \left(\sum_{a=1}^{h_{\bm{r}}}L_d(a)^q\right)^{1/q}\\
        \label{eq:group-prefix-substitution}
        &\leq 2^{d-1}\left(\prod_{i=1}^dK_i\right)
              2^{\gamma(R+d-1)}h_{\bm{r}}^\gamma\\
        \label{eq:group-cost}
        &=2^{(d-1)(1+\gamma)}
              \left(\prod_{i=1}^dK_i\right)v_{\bm{r}}^{\gamma}.
    \end{align}
    The factor $2^{(d-1)(1+\gamma)}$ combines $2^{d-1}$ from the line-weight bound with $2^{\gamma(d-1)}$ from the $d-1$ dyadic prefix cutoffs, using $v_{\bm{r}}=2^Rh_{\bm{r}}$.
    
    \paragraph{Pack disjoint boxes.}
    We bound the sum of the $v_{\bm{r}}$ by comparing them with the sizes of disjoint boxes inside $A$.
    For each $\bm{r}\in\mathcal R_A$, choose $\bm{\mu}\in S_{\bm{r}}$ with $\mu_d+1=h_{\bm{r}}$.
    Since $\bm{\mu}\in S_{\bm{r}}$, we have $\mu_i+1\in E_{r_i}$, hence $\mu_i\geq2^{r_i}-1$ for $1\leq i\leq d-1$.
    Lemma~\ref{lem:line} also gives $\bm{\mu}\in A$.
    Set $T_0=\{0\}$ and $T_r=\{2^{r-1},\ldots,2^r-1\}$ for $r\geq1$.
    Since $\max T_{r_i}=2^{r_i}-1\leq\mu_i$ for $1\leq i\leq d-1$ and $h_{\bm{r}}-1=\mu_d$, every index in the following box lies coordinatewise below $\bm{\mu}$.
    Downward closure therefore gives
    \begin{equation}
        \label{eq:retained-box}
        Q_{\bm{r}}=
              \left(\prod_{i=1}^{d-1}T_{r_i}\right)
              \times\{0,\ldots,h_{\bm{r}}-1\}\subseteq A.
    \end{equation}
    The intervals $T_0,T_1,T_2,\ldots$ occupy disjoint index ranges.
    Different groups differ in at least one of their first $d-1$ ranges, so their retained boxes are disjoint even if the heights overlap.
    Since $|T_r|=2^{r-1}$ for $r\geq1$ and $|T_0|=1$, we have $2^r\leq2|T_r|$ for every $r\in\Nzero$.
    Thus, for every $\bm{r}\in\mathcal R_A$,
    \begin{equation}
        \label{eq:group-volume-comparison}
        v_{\bm{r}}
        =h_{\bm{r}}\prod_{i=1}^{d-1}2^{r_i}
        \leq2^{d-1}h_{\bm{r}}\prod_{i=1}^{d-1}|T_{r_i}|
        =2^{d-1}|Q_{\bm{r}}|.
    \end{equation}
    The boxes are pairwise disjoint subsets of $A$, so summing this comparison gives
    \begin{equation}
        \label{eq:packing}
        \sum_{\bm{r}\in\mathcal R_A}v_{\bm{r}}
        \leq2^{d-1}\sum_{\bm{r}\in\mathcal R_A}|Q_{\bm{r}}|
        =2^{d-1}\left|\bigcup_{\bm{r}\in\mathcal R_A}Q_{\bm{r}}\right|
        \leq2^{d-1}|A|=2^{d-1}N.
    \end{equation}

    \paragraph{Sum the group bounds.}
    Finally, we use \eqref{eq:absolute-cost}, \eqref{eq:group-cost} and \eqref{eq:packing}.
    Since $\gamma\geq1$, this gives
    \begin{align}
        \label{eq:group-sum}
        \Lambda(P_A)
        &\leq 2^{(d-1)(1+\gamma)}\left(\prod_{i=1}^dK_i\right)
                           \sum_{\bm{r}\in\mathcal R_A}v_{\bm{r}}^{\gamma}\\
        \label{eq:finish-transfer}
        &\leq 2^{(d-1)(1+\gamma)}\left(\prod_{i=1}^dK_i\right)
                           \left(\sum_{\bm{r}\in\mathcal R_A}v_{\bm{r}}\right)^\gamma
        \leq 2^{(d-1)(1+2\gamma)}\left(\prod_{i=1}^dK_i\right)N^\gamma.
    \end{align}
    This proves the assertion.
\end{proof}

By Lemma~\ref{lem:moment}, the assumptions of Proposition~\ref{prop:transfer} hold in each coordinate with $q=q_d$, $\gamma=\alpha_{q_d}/q_d=\betaL$, and $K_i=K_{q_d}^{1/q_d}$. 
This finally proves Theorem~\ref{thm:main}.

\section{Numerical experiments}
\label{sec:numerical-validation}

We now compare the upper bounds from Theorem~\ref{thm:main} with numerical estimates of the Lebesgue constants for four families of lower sets in dimensions two and three, namely total-degree sets $\sum_{i=1}^d \nu_i\leq n$, weighted-degree sets $\sum_{i=1}^d i \nu_i \leq n$, hyperbolic crosses $\prod_{i=1}^d(\nu_i+1)\leq T$, and tensor cubes $\{0,\ldots,k-1\}^d$.
All grids use, in every coordinate, the $\mathfrak R$-Leja sequence obtained from the canonical bit-reversal disk Leja sequence starting at $1$ \cite[equations~(2.2)--(2.3)]{Chkifa2015Bounds}.
Their sizes are chosen to approximate the target node counts
\begin{equation}
    \label{eq:numerical-targets}
    B_j=\left\lfloor 2^{8+j/2}+\tfrac12\right\rfloor,
    \qquad 
    j\in\{0,\ldots,14\},
\end{equation}
which increase by a factor $\sqrt2$ before rounding, from $256$ to $32\,768$.
For each family and target $B_j$, the selected grid minimizes $|\log_2(N/B_j)|$ subject to $256\leq N\leq32\,768$, and the plots use its actual node count $N$.

For each non-tensor grid, evaluations of the Lebesgue function at $16\,384$ scrambled Sobol points and $2\,048$ boundary points are followed by six bounded local Powell searches.
Tensor cubes admit the factorization $\Lambda(P_A)=L(k)^d$, which we evaluate numerically for every integer prefix length $k$ with $256\leq k^d\leq32,768$, so that the plotted curves resolve the oscillations between successive prefixes.
In both cases, the computed values $\widehat\Lambda_A$ are numerical estimates from below, subject to roundoff.

\begin{figure}[!htbp]
    \centering
    \vspace{0.5em}
    \includegraphics[width=\linewidth]{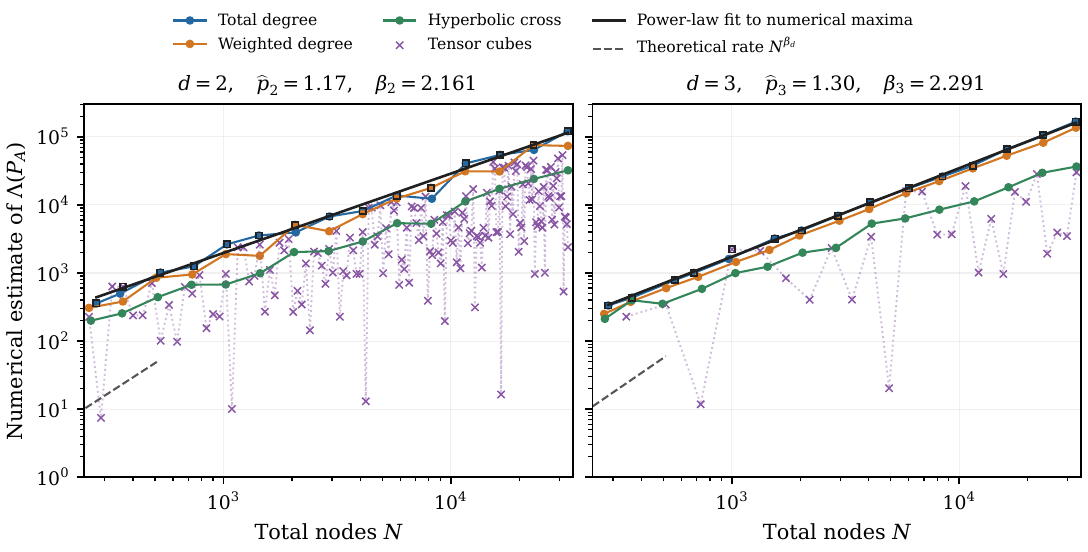}
    \caption{
        Numerical growth of the Lebesgue constants for \(256 \leq N \leq 32\,768\). 
        The tensor-product curves include all possible integer prefix lengths. 
        Black squares mark the largest estimate among the four grids matched to each of the fifteen target node counts, while solid black lines show the corresponding power-law fits. 
        Dashed lines indicate the theoretical slopes \(\beta_d\), their vertical placement is arbitrary.
    }
    \label{fig:numerical-growth}
\end{figure}

Figure~\ref{fig:numerical-growth} and Table~\ref{tab:numerical-slopes} summarize the results.
We plot $\log\widehat\Lambda_A$ against $\log N$ for each family, using the fifteen matched grids for each non-tensor family and every prefix for tensor cubes.
For the pooled plot, we select the largest estimate among the four grids matched to each $B_j$ and retain that grid's actual $N$.

\begin{table}[H]
    \centering
    \begin{tabular}{@{}lrr@{}}
        \toprule
        Family & $d=2$ & $d=3$\\
        \midrule
        Total degree & $1.159$ & $1.310$\\
        Weighted degree & $1.164$ & $1.309$\\
        Hyperbolic cross & $1.056$ & $1.078$\\
        Tensor cube & $0.956$ & $0.927$\\
        \midrule
        Pooled maxima & $1.166$ & $1.301$\\
        \bottomrule
    \end{tabular}
    \caption{
        Finite-range exponents obtained from linear fits of $\log\widehat\Lambda_A$ against $\log N$, rounded to three decimal places.
        The fits for tensor cubes use $166$ data pairs for $d=2$ and $26$ for $d=3$, while all other fits use fifteen.
        The last row reports the exponents of the pooled fits shown in Figure~\ref{fig:numerical-growth}.
    }
    \label{tab:numerical-slopes}
\end{table}

Over the tested range, the fitted exponents vary across the four families but are all substantially smaller than the theoretical exponent $\beta_d$ in Theorem~\ref{thm:main}, as shown in Figure~\ref{fig:numerical-growth} and Table~\ref{tab:numerical-slopes}.

\section{Conclusion}
\label{sec:conclusion}

Using the known binary estimates for $\mathfrak R$-Leja prefixes, we established the univariate prefix-moment bound~\eqref{eq:moment} in Lemma~\ref{lem:moment}.
We then bounded the total absolute coefficient weight along each coordinate line in Lemma~\ref{lem:line} and the weighted products in Lemma~\ref{lem:holder}.
Together with~\eqref{eq:absolute-cost} and~\eqref{eq:packing}, these estimates yielded Proposition~\ref{prop:transfer}, which transfers an algebraic bound on the $q$th root of the univariate prefix-moments to an algebraic bound on the multivariate Lebesgue constant with the same exponent in the total node count.
Applying this proposition with $q=d/(d-1)$ and $\gamma=\alpha_q/q=\beta_d$, we proved Theorem~\ref{thm:main}, i.e., for each fixed dimension $d\geq2$ and every nonempty finite lower set $A\subset\Nzero^d$, we have $\Lambda(P_A)\leq C_dN^{\beta_d}$, where $N=|A|$ and $C_d$ depends only on $d$.
The exponent we obtained in~\eqref{eq:beta} takes the values $\beta_2=2.160964\ldots$ in two dimensions and $\beta_3=2.291167\ldots$ in three dimensions.

With this result, we improved the cubic estimate~\eqref{eq:cubic} and reduced the exponent $3-1/d$ obtained from quadratic prefix bounds in~\eqref{eq:quadratic-md}.
In particular, we obtained a sharper bound for total-degree sets compared with~\eqref{eq:shape-total-degree-bound}.
For dyadic sparse sets, the direct box estimate~\eqref{eq:shape-dyadic-bound} combined with the sharper dyadic-prefix bounds yields the bound $O_d(N)$.
By contrast, our main theorem applies to \emph{arbitrary lower sets} without any restrictions on their shape.

Over the tested range, the fitted exponents of the numerically computed Lebesgue-constant estimates in Section~\ref{sec:numerical-validation} were substantially smaller than the theoretical exponent $\beta_d$ for all four families.
These observations motivate seeking sharper bounds. 
In particular, any reduction of the exponent in the univariate prefix-moment bound~\eqref{eq:moment} with $q=q_d$ would, using Proposition~\ref{prop:transfer}, yield a smaller exponent in the uniform multivariate bound of Theorem~\ref{thm:main}.

\section*{Acknowledgements}

This work was partially funded by the Center for Advanced Systems Understanding (CASUS), financed by Germany's Federal Ministry of Research, Technology and Space (BMFTR) and by the Saxon Ministry for Science, Culture and Tourism (SMWK) with tax funds on the basis of the budget approved by the Saxon State Parliament.

\begingroup
\setstretch{1.0}
\bibliographystyle{plainnat}
\bibliography{references}

@article{CarnicerGodes2015,
  author = {Carnicer, Jes{\'u}s and God{\'e}s, Carmen},
  title = {Lagrange polynomials of lower sets},
  journal = {Dolomites Research Notes on Approximation},
  volume = {8},
  number = {Special Issue},
  pages = {1--10},
  year = {2015},
  doi = {10.14658/PUPJ-DRNA-2015-Special_Issue-1}
}

@incollection{Chkifa2015Bounds,
  author = {Chkifa, Abdellah},
  title = {New Bounds on the {Lebesgue} Constants of {Leja} Sequences on the Unit Disc and on {$\mathfrak R$}-{Leja} Sequences},
  booktitle = {Curves and Surfaces},
  series = {Lecture Notes in Computer Science},
  volume = {9213},
  pages = {109--128},
  year = {2015},
  publisher = {Springer},
  doi = {10.1007/978-3-319-22804-4_9},
  url = {https://arxiv.org/abs/1503.01731v2},
  note = {Numbered locators refer to the arXiv version, 16 May 2015}
}

@article{ChkifaCohenSchwab2014,
  author = {Chkifa, Abdellah and Cohen, Albert and Schwab, Christoph},
  title = {High-dimensional adaptive sparse polynomial interpolation and applications to parametric {PDE}s},
  journal = {Foundations of Computational Mathematics},
  volume = {14},
  number = {4},
  pages = {601--633},
  year = {2014},
  doi = {10.1007/s10208-013-9154-z}
}

@article{Chkifa2013Lebesgue,
  author = {Chkifa, Moulay Abdellah},
  title = {On the {Lebesgue} constant of {Leja} sequences for the complex unit disk and of their real projection},
  journal = {Journal of Approximation Theory},
  volume = {166},
  pages = {176--200},
  year = {2013},
  doi = {10.1016/j.jat.2012.11.005}
}

@article{Chkifa2022Newton,
  author = {Chkifa, Moulay Abdellah},
  title = {Newton interpolation using {$\Re$}-{Leja} sequences},
  journal = {BIT Numerical Mathematics},
  volume = {62},
  pages = {1461--1485},
  year = {2022},
  doi = {10.1007/s10543-022-00918-y}
}

@incollection{CohenChkifa2015,
  author    = {Cohen, Albert and Chkifa, Abdellah},
  title     = {On the Stability of Polynomial Interpolation Using Hierarchical Sampling},
  booktitle = {Sampling Theory, a Renaissance},
  editor    = {Pfander, G{\"o}tz E.},
  series    = {Applied and Numerical Harmonic Analysis},
  pages     = {437--458},
  publisher = {Birkh{\"a}user},
  year      = {2015},
  doi       = {10.1007/978-3-319-19749-4_12},
  url       = {https://hal.science/hal-01353241v1},
}

@article{DynFloater2014,
  author = {Dyn, Nira and Floater, Michael S.},
  title = {Multivariate polynomial interpolation on lower sets},
  journal = {Journal of Approximation Theory},
  volume = {177},
  pages = {34--42},
  year = {2014},
  doi = {10.1016/j.jat.2013.09.008}
}

@article{KangWilcox2017,
  author = {Kang, Wei and Wilcox, Lucas C.},
  title = {Mitigating the curse of dimensionality: sparse grid characteristics method for optimal feedback control and {HJB} equations},
  journal = {Computational Optimization and Applications},
  volume = {68},
  number = {2},
  pages = {289--315},
  year = {2017},
  doi = {10.1007/s10589-017-9910-0},
  url = {https://arxiv.org/abs/1507.04769}
}

@article{MininiWasem2025,
  author = {Minini, Jocelyn and Wasem, Micha},
  title = {Counting nodes in {Smolyak} grids},
  journal = {Ars Combinatoria},
  volume = {162},
  pages = {149--157},
  year = {2025},
  doi = {10.61091/ars162-11}
}

@article{SickelUllrich2007,
  author = {Sickel, Winfried and Ullrich, Tino},
  title = {The {Smolyak} algorithm, sampling on sparse grids and function spaces of dominating mixed smoothness},
  journal = {East Journal on Approximations},
  volume = {13},
  number = {4},
  pages = {387--425},
  year = {2007},
  url = {https://www.tu-chemnitz.de/mathematik/ang_analysis/Publications/TW07.pdf}
}

@article{
    Hecht2026,
    author = {Hecht, Michael and Hofmann, Phil-Alexander and Wicaksono, Damar and Hernandez Acosta, Uwe and Gonciarz, Krzysztof and Kissinger, Jannik and Sivkin, Vladimir and Sbalzarini, Ivo F},
    title = {Multivariate {N}ewton interpolation in downward closed spaces reaches the optimal {B}ernstein–{W}alsh approximation rate},
    journal = {IMA Journal of Numerical Analysis},
    pages = {draf137},
    year = {2026},
    month = {03},
    issn = {0272-4979},
    doi = {10.1093/imanum/draf137},
    url = {https://doi.org/10.1093/imanum/draf137},
    eprint = {https://academic.oup.com/imajna/advance-article-pdf/doi/10.1093/imanum/draf137/67244017/draf137.pdf},
}

@article{CalviPhung2011,
  author  = {Calvi, Jean-Paul and Phung, Van Manh},
  title   = {On the {Lebesgue} constant of {Leja} sequences for the unit disk and its applications to multivariate interpolation},
  journal = {Journal of Approximation Theory},
  volume  = {163},
  number  = {5},
  pages   = {608--622},
  year    = {2011},
  doi     = {10.1016/j.jat.2011.02.001}
}

@article{CalviPhung2012,
  author  = {Calvi, Jean-Paul and Phung, Van Manh},
  title   = {Lagrange interpolation at real projections of {Leja} sequences for the unit disk},
  journal = {Proceedings of the American Mathematical Society},
  volume  = {140},
  number  = {12},
  pages   = {4271--4284},
  year    = {2012},
  doi     = {10.1090/S0002-9939-2012-11291-2}
}
\endgroup
\end{document}